\documentclass[12pt,a4paper]{article}
\usepackage{amsfonts}
\usepackage{graphicx}

\newfont{\bbf}{msbm10 at 12pt}
\newfont\bbfsm{msbm10 at 9pt}

\def\N{\mathbb {N}}

\def\Z{\mathbb {Z}}
\def\R{\mathbb {R}}
\def\C{\mathbb {C}}
\def\Cbar{\overline{\C}}
\def\Cstar{\C^*}

\newfont{\script}{eusm10 at 12pt}
\newfont{\scriptsmall}{eusm10 at 9pt}

\def\ovl{\overline}
\def\phi1{\phi}
\def\phi{\varphi}
\def\eps{\varepsilon}
\def\theta{\vartheta}
\def\Re{\mbox{\rm Re}}
\def\Im{\mbox{\rm Im}}

\def\sm{\setminus}

\def\*{ {\mbox{\tt $\star$}} }
\def\0{ {\mbox{\tt 0}} }
\def\1{ {\mbox{\tt 1}} }
\def\2{ {\mbox{\tt 2}} }
\def\3{ {\mbox{\tt 3}} }
\def\4{ {\mbox{\tt 4}} }

\def\St#1#2 {{\small\tt\rule{0pt}{0pt}_{\mbox{#1}}^{\mbox{#2}} }}

\def\<{\prec}
\def\>{\succ}

\ifx\ThmNum\undefined 
   \newtheorem{theorem}{Theorem}[section]
   
\else 
   \newtheorem{theorem}{Theorem} 
   
\fi

\newtheorem{proposition}[theorem]{Proposition}
\newtheorem{lemma}[theorem]{Lemma}
\newtheorem{definition}[theorem]{Definition}

\newtheorem{corollary}[theorem]{Corollary}

\def\proof{\par\medskip\noindent {\sc Proof. }}
\def\proofof #1 {\par\medskip\noindent {\sc Proof of #1. }}
\def\sketch{\par\medskip\noindent{\sc Sketch of proof. }}
\def\sketchof #1 {\par\medskip\noindent {\sc Sketch of proof of #1. }}
\def\Box{\framebox[10pt]{\rule{0pt}{3pt}}}
\def\nix{\rule{0pt}{2pt}}
\def\qed{\qedd\par\medskip\noindent}
\def\qedd{\nix\nolinebreak\hfill\hfill\nolinebreak$\Box$}
\def\remark{\par\medskip \noindent {\sc Remark. }}
\def\lineclear
    {\rule{0pt}{0pt}\nopagebreak\par\nopagebreak\noindent}

\newlength\captionwidth
\def\LabelCaption#1#2{
   \centerline{\parbox{\captionwidth}{   
   \caption{\sl #2}  \label{#1} } }
}

\def\reminder #1 {{\sf #1}}
\def\hide #1 {}
\long\def\longhide #1 {}

\def\A{\mathcal{A}}
\def\Sym{\mathcal{S}}
\def\su{{\underline{s}}}
\def\sigsu{{\sigma(\su)}}
\def\overview{\par\medskip\noindent {\sc Overview: }}
\newcommand\Intro[1] {{\small #1}\par\medskip}
\def\Re{{\rm Re}}
\def\Im{{\rm Im}}
\def\Cp{{\C\,'}}
\begin{document}
\title{Escaping Points of the Cosine Family}
\author{G\"unter Rottenfu\ss er and Dierk Schleicher}


\maketitle

\begin{abstract}
We study the dynamics of iterated cosine maps $E\colon z \mapsto ae^z+be^{-z},$
with $a,b \in \C \sm \{0\}$. We show that the points which converge to
$\infty$ under iteration are organized in the form of rays and, as 
in the exponential family, every escaping point is either on one of these
rays or the landing point of a unique ray. Thus we get a complete classification
of the escaping points of the cosine family, confirming a conjecture of
Eremenko in this case. We also get a particularly strong version of the
``dimension paradox'': the set of rays has Hausdorff dimension $1$, while the
set of points these rays land at has not only Hausdorff dimension $2$ but
infinite planar Lebesgue measure.
\end{abstract}

\tableofcontents

\section{Introduction}
\label{SecIntro}

The dynamics of iterated polynomials has been investigated quite successfully,
particularly in the past two decades. The study begins with a description of
the {\em escaping points:} those points which converge to $\infty$ under
iteration. It is well known that the set of escaping points is an open
neighborhood of $\infty$ which can be parametrized by {\em dynamic rays}. The
Julia set can then be studied in terms of landing properties of dynamic rays.

For entire transcendental functions, the point $\infty$ is an essential
singularity (rather than a superattracting fixed point as for polynomials).
This makes the investigation of the dynamics much more difficult. In
particular, there is no obvious structure of the set of escaping points. 
Eremenko~\cite{E} showed that for every entire transcendental function, the
set of escaping points is always non-empty, and he asked whether each
connected component (or even every path component) was unbounded. In
\cite{Escaping}, this question was answered in the affirmative for the special
case of exponential functions $z\mapsto\lambda e^z$: every escaping point can
be connected to $\infty$ along a unique curve running entirely along
escaping points. In this paper, we extend this result to the family of
cosine maps $E_{a,b}\colon z\mapsto ae^z+be^{-z}$ with $a,b\in\Cstar$. 

In many ways, the escaping points of the cosine family behave quite similarly
to those of the exponential family; therefore, the present paper is very
similar to \cite{Escaping}. It is based on the Diploma thesis \cite{Guenter}.
Based on experience with these and other families of maps, we believe that
similar results should hold for much larger classes of entire transcendental
functions, possibly of bounded type (which means that the set of asymptotic or
critical values is bounded). This paper is thus a contribution to the program
to make polynomials tools, in particular dynamic rays, available for the study
of iterated entire transcendental functions. 

Our main result is a classification of escaping points for every map
$E_{a,b}$; this classification is the same for all such maps (with
natural exceptions if one or both critical orbits escape). A byproduct is an
affirmative answer to Eremenko's question as mentioned above: every path 
connected component of the set of escaping points is a curve starting at
$\infty$. Similarly as for exponential functions, but quite unlike the
polynomial case, certain of these curves land at points in $\C$ which are
escaping points themselves. A {\em dynamic ray} is a connected component of
the escaping set, removing the landing points (for those curves which land at
escaping points). It turns out that the union of all the uncountably many
dynamic rays has Hausdorff dimension $1$ (in analogy to results
\cite{Karpinska,Escaping} for the exponential family). However, by a result of
McMullen~\cite{McMullen}, the set of escaping points in the cosine family has
infinite planar Lebesgue measure (this is one main difference to the
exponential case, where the set of escaping points has zero measure).
Therefore, the entire measure of the escaping set sits in the landing points
of those rays which land at escaping points.

In \cite{Cosine}, these results are extended even further: if both critical
orbits of $E_{a,b}$ are strictly preperiodic, then the set $R$ of dynamic rays
still has Hausdorff dimension $1$ and every dynamic ray lands somewhere in
$\C$. Therefore, $\C\sm R$ is ``most of $\C$'' (the complement of a
one-dimensional set). It turns out that each $z\in\C\sm R$ is the landing point
of one or several of the rays in $R$: the rays in the one-dimensional set $R$
manage to connect all the remaining points to $\infty$ by curves in $R$! This
highlights another difference between the dynamics in the cosine family and
in the exponential family: in the exponential family, there is no case known
where the Julia set is the entire complex plane and every dynamic ray lands
\cite[Theorem~1.6]{Lasse}; this is due to the enormous contraction in the
asymptotic tract of the asymptotic value $0$, while the singular values in the
cosine family are simply two critical values with much better-behaved
properties.

We start this paper in Section~\ref{SecEscaping} by setting up a partition for 
symbolic dynamics. In Section~\ref{SecRayTails}, we construct {\em ray tails}
which are curves of escaping points terminating at $\infty$. These ray tails
are extended to entire {\em dynamic rays} in Section~\ref{SecRays}. In
Section~\ref{SecHorizontalEscape}, we prove that points on the same dynamic ray
move away from each other very quickly, and this leads to a complete
classification of escaping points in Section~\ref{SecClassification}. Finally,
in Section~\ref{SecEpilogue} we discuss the implications in terms of Hausdorff
dimension.

{\sc Acknowledgment}. 
We would like to thank Alexandra Kaffl and Johannes R\"uckert and in
particular Markus F\"orster and Lasse Rempe for helpful comments on earlier
drafts of this paper. We are grateful to the Institut Henri Poincar\'e in
Paris for its hospitality while this paper was finished, and to Adrien Douady
and Hans Henrik Rugh for having organized the program which brought us there
in the fall of 2003. The first author would like to thank the DAAD for its
support during the stay.

{\sc Notation.}
We consider the maps 
\[
E_{a,b}(z):=ae^z+be^{-z} \mbox{ for } a,b \in \Cstar:=\C \sm \{0\}
\]
and their iterates 
$E_{a,b}^{\circ n}(z)$. Usually we will omit the parameters $a$ and $b$
and write $E(z)$ for $E_{a,b}(z)$.
Set $c:=\frac{1}{2}\ln\left(\frac{b}{a}\right)$, where the branch of the
logarithm is chosen such that $|\Im(c)| \leq \pi  /2$.
The critical points of $E$ are
\[
C_{crit}=\left\{ c+i\pi n, n\in  \Z\right\}
\]
and the critical values $v_{1/2}= \pm 2\sqrt{ab}$, choosing signs so that
$v_1$ is the image of $c+2\pi i\Z$, while $v_2$ is the image of $c+i\pi+2\pi
i\Z$. There are no asymptotic values in $\C$. We will use the notation $\Cstar
:= \C \sm \{0\}$ and $\Cbar:=\C \cup \{\infty\}$. Also, we will need
$F \colon \R^+_0 \rightarrow \R^+_0,\,\,  F(t):=\exp(t)-1$.
Now let $\alpha:=\ln a$ and $\beta:=\ln b$, choosing branches so that
$|\Im(\alpha)|\leq \pi$, $|\Im(\beta)|\leq\pi$. Let furthermore $K:=\min \left\{ |a|,|b| \right\}$ , $K_{\max}:=\max \left\{ |a|,|b| \right\}$ and $M:=\max\{|\alpha|,|\beta|\}$.
Given $a,b\in\Cstar$, let
\begin{eqnarray}
T_{a,b} &\geq&
\max\left\{\rule{0pt}{24pt}
\sqrt{\left|\frac{2b}{a}\right|}(|a|+|b|),
\sqrt{\left|\frac{2a}{b}\right|}(|a|+|b|), 8  |ab|, 
1,\right.
\nonumber\\
&&\qquad\quad\left. 
\frac 1 2 \ln\left|\frac {2b} a \right|,\frac 1 2 \ln \left|\frac {2a}
b\right|,
\ln \frac{4}{|a|}, \ln \frac{4}{|b|}
\rule{0pt}{24pt}\right\}
 +M+2 
\label{EqDefTab}
\end{eqnarray}
be the least value for which $F(T_{a,b})\geq T_{a,b}+M+4$. Note that 
\begin{equation}
2\sqrt{|ab|} \leq \max\{8|ab|,1\} \,\,:
\label{EqInequTab}
\end{equation}
if $|ab|<1/16$, then $2\sqrt{|ab|}<1/2<1$, while if $|ab|\geq 1/16$, then
$2\sqrt{|ab|}\leq 8|ab|$. Thus $T_{a,b}\geq 2\sqrt{|ab|}+M+2$.

We will use the following sets (their significance will be explained in
Section~\ref{SecEscaping} and Figure~\ref{FigPartition}): if $\Im(v_1)\geq
\Im(v_2)$, then set
\begin{eqnarray}
\A &:=& \left\{ \strut z\in \C \colon 
z=\lambda v_1+(1-\lambda) v_2; \lambda\in[0,1]\right\} \nonumber\\ 
&&\cup
\left\{\strut z\in \C \colon \Re(z)=\Re(v_1), \Im(z) \geq \Im(v_1) \right\}
\,\,;
\label{EqDefA}
\end{eqnarray}
this is the segment between $v_1$ and $v_2$, together with the vertical ray
starting at $v_1$ in upwards direction. If $\Im(v_1)<\Im(v_2)$, we reverse
the last inequality in the definition of $\A$: the last ray is replaced by the
downwards vertical ray at $v_1$. In both cases, set $\Cp := \C \sm \A$.

\section{Escaping Points and Symbolic Dynamics}
\label{SecEscaping}

\Intro{In this section, we set up a partition of the complex plane and
define symbolic dynamics of escaping points.}

\begin{definition}[Escaping Points]
\label{escpts} \lineclear
A point $z\in \C$ with $|E^{\circ k}(z)| \to \infty$ for $k
\to \infty$ is called an {\em escaping point} and its orbit is called {\em
   escaping orbit}.
\end{definition}

\begin{lemma}[Real Parts of Escaping Orbits]
\label{realescorb}\lineclear
If $(z_k)$ is an orbit with  $|z_k|\to\infty$ for $k\to\infty$,
then $|\Re(z_k)| \to \infty$.
\end{lemma}
\proof
This follows from the standard estimate
\begin{equation}
|z_{k+1}|\leq|a|\exp(\Re(z_k))+|b|\exp(\Re(-z_k)) \,\,	.
\label{EqStandardEstimate}
\end{equation}\qed

In order to introduce symbolic dynamics we need a useful partition of
the complex plane. First define
\[
\Z_S:=\Z \times \{L,R\},~
\Z_R:=\Z \times \{R\} \mbox{ and }
\Z_L:=\Z \times \{L\},
\]
and we write $n_R:=(n,R)$ and $n_L:=(n,L)$ for $n \in \Z$.
Let us first consider the ``straight partition''
(Figure~\ref{FigFirstPartition}): for $j\in \Z_R$ define the strip $R'_j$ as
follows
  \[
R'_j:=\left\{z \in \C\ \colon \Im(c)+2\pi j < \Im(z) <
  \Im(c)+2\pi(j+1); \Re(c)<\Re(z) \right\}
\]
  and for $j\in \Z_L$
\[
R'_j:=\left\{z\in \C \colon \Im(c)+2\pi j < \Im(z)< \Im(c)+ 2\pi(j+1);
    \Re(z)<\Re(c)\right\}.
\]

Then the $R'_j$ form a rather simple partition of $\C$. The image of the
boundary of each strip is the part of the straight line through the critical
values ending in $v_1$. It will be convenient to modify the partition 
so that the
real part of the image of the boundary is bounded. To introduce our partition
(see Figure~\ref{FigPartition}) now define $\A$ as in (\ref{EqDefA}) and set
$\Cp := \C \sm \A$. Then define the strips $R_j$ as connected components of
$E^{-1}(\C \sm \A)$, so that
\[
E\colon R_j \rightarrow \Cp
\]
is a conformal isomorphism for all $ j \in \Z_S$. For $s\in \Z_S$ we
denote the inverse mapping of
$E\colon R_s \rightarrow \Cp$ by $L_s \colon \Cp \rightarrow R_s$.
Note that the strips are open and that the union
of their closures is $\C$. Label the strips such that $\ovl R_{0_R}$ contains
some right end of $\R^+$, $\ovl R_{0_L}$ contains some left end of $\R^-$ ,
and that $z\in R_{(n,L)}$ iff $z+2\pi i\in R_{(n+1,L)}$ and $z\in R_{(n,R)}$
iff $z+2\pi i\in R_{(n+1,R)}$. The following property of the strips will be
used throughout:

\begin{lemma}[Height of Strips Bounded by $3\pi$]
\label{LemHeightStrips} \lineclear
If $z, w \in R_{n,L}\cup R_{n,R}$, then $|\Im(z)-\Im(w)|<3\pi$ and
$|\Im(z)-2\pi i n|<3\pi$.
\end{lemma}

\proof
The straight line through $v_1$ and $v_2$ divides $\C$ into two halfplanes,
and the preimages under $E_{a,b}$ of each of them are horizontal half-strips
of height $\pi$. Therefore, each $R_s$ is contained in three of those
half-strips. Since $\ovl R_{0_L}$ and $\ovl R_{0_R}$ intersect the real axis, 
the second part follows immediately.
\qed

\begin{figure}
\includegraphics[width=0.9\textwidth]{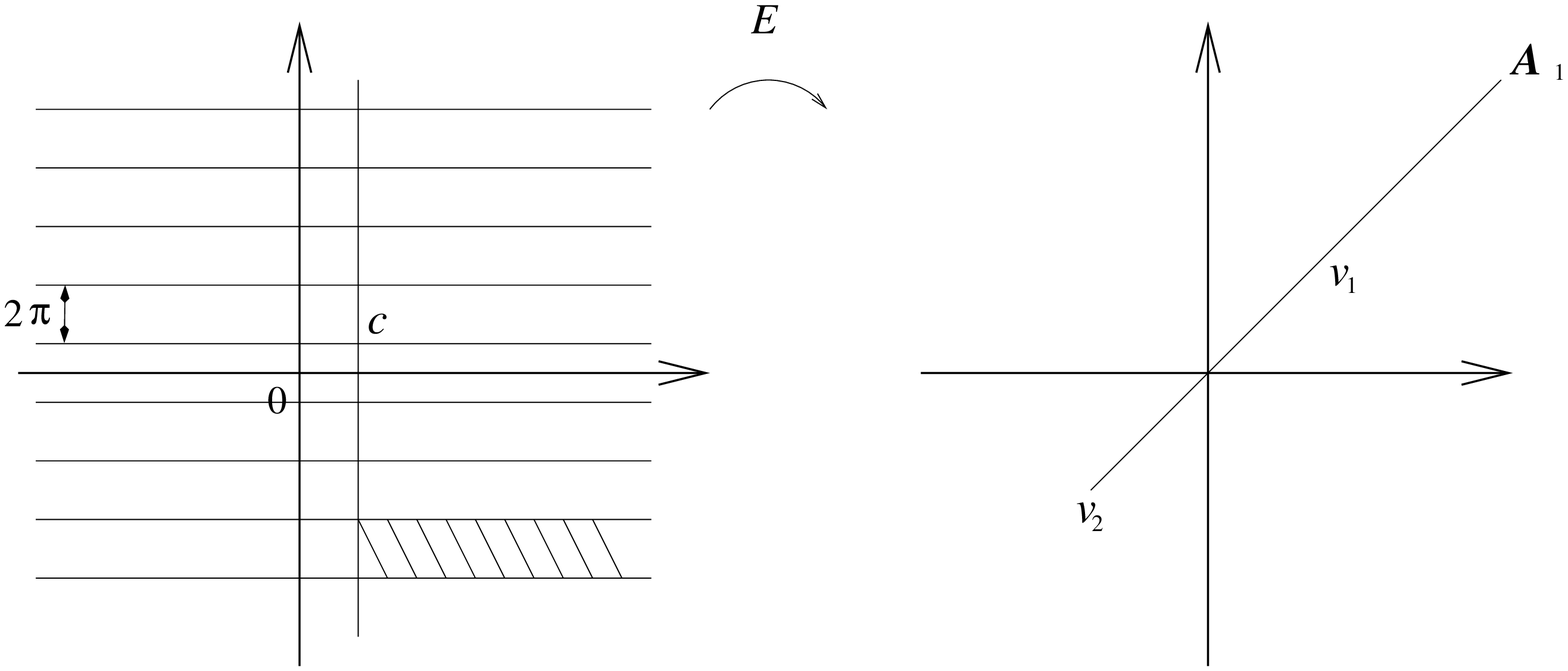}
\LabelCaption{FigFirstPartition}{The partition formed by the $R'_j$: the strips
have simple shapes, but the real part of $\A_1$ is unbounded.}

\includegraphics[width=0.9\textwidth]{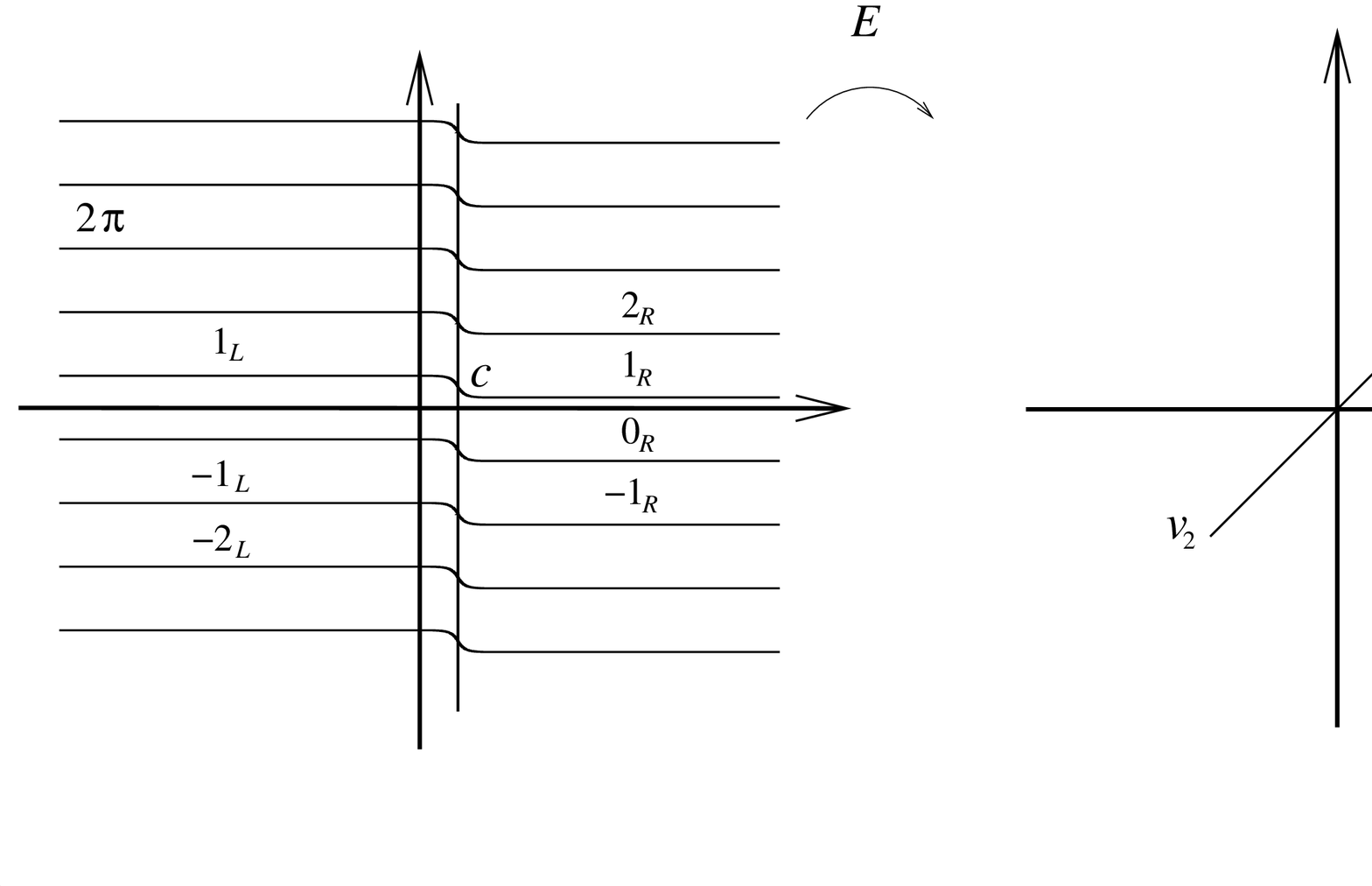}
\LabelCaption{FigPartition}{The partition we use: the set $\A$ has bounded real
parts.}

\includegraphics[width=0.9\textwidth]{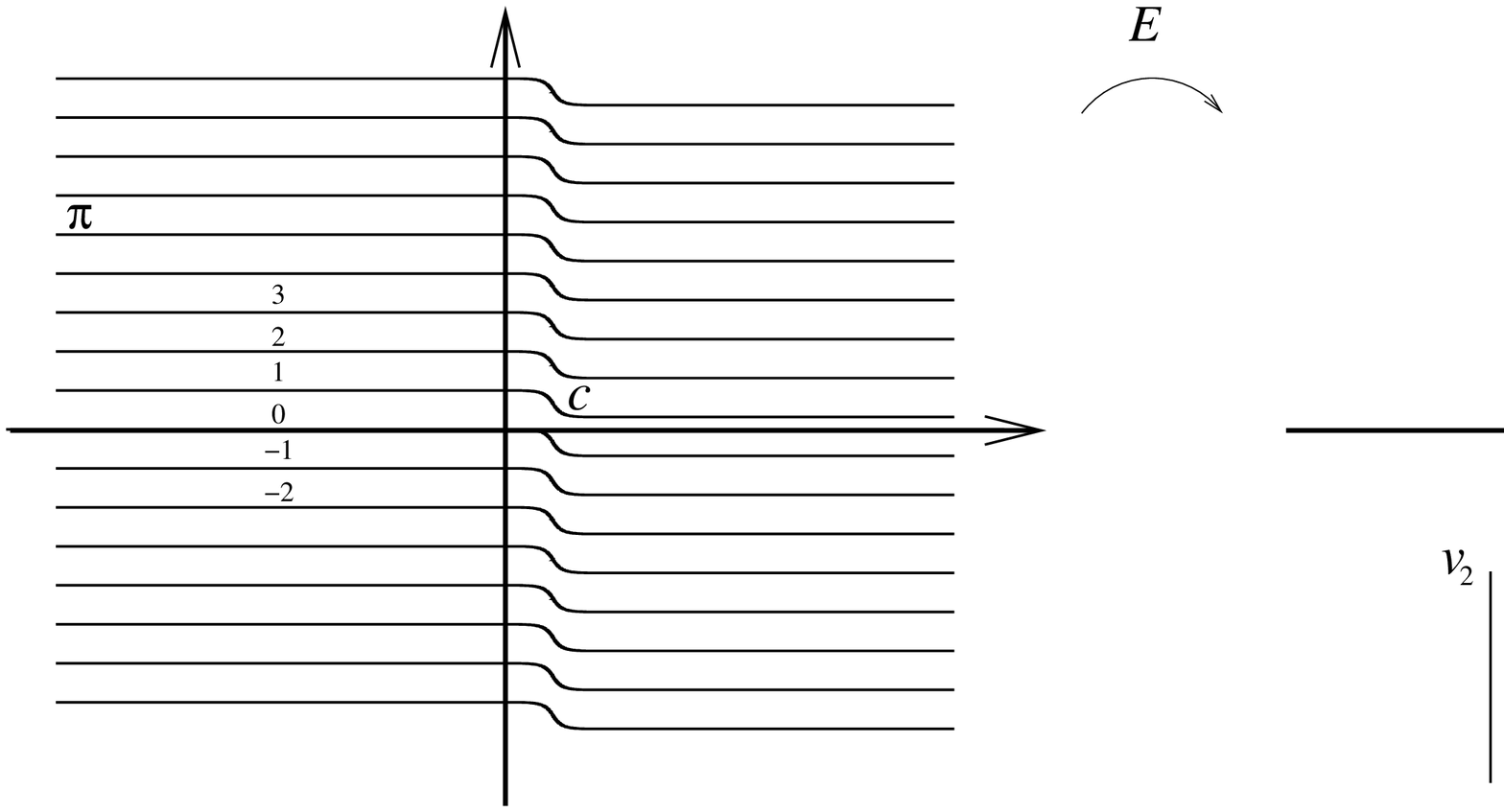}
\LabelCaption{FigSecondPartition}{A partition with bounded real parts of $\A_2$
which could be used.}
\end{figure}

\begin{lemma}[Sign of Real Part]
\label{LemSignRealPart} \lineclear
If $|E(z)|>|a|+|b|$ for some $z\in\C$, then $\Re(z)>0$ iff $z\in R_{n,R}$ for
some $n\in\Z$, and $\Re(z)<0$ iff $z\in R_{n,L}$ for some $n\in\Z$.
\end{lemma}
\proof
Suppose that $\Re(c)<0$. Then $\{z\in R_{n,L}\colon \Re(z)>0\}$ is empty for
all $n$, while every $z\in R_{n,R}$ with $\Re(z)<0$ has real parts between
$\Re(c)$ and $0$. The vertical line through $c$ (containing all the critical
points) is mapped under $E$ to the segment connecting the two critical values
$\pm2\sqrt{ab}$, while the imaginary axis maps to an ellipse with major axis
$|a|+|b|$. Therefore if $|E(z)|>|a|+|b|$, then $|\Re(z)|$ has the same
sign as the unbounded part of the strip $R_{(n,L)}$ or $R_{(n,R})$ containing
$z$.

The case $\Re(c)\geq 0$ is analogous.
\qed

There is a number of further conceivable partitions, such as the one in
Figure~\ref{FigSecondPartition}, all with their particular advantages
but with a different syntax of symbolic sequences.

\begin{definition}[External Address]
\label{extadre} \lineclear
Let $\Sym:=\Z_S^\N=\{(s_1s_2s_3 \ldots)\colon s_k \in \Z_S \}$
be the sequence space over $\Z_S$ and let $\sigma \colon \Sym \rightarrow
\Sym, (s_1s_2s_3s_4 \ldots) \mapsto (s_2s_3s_4 \ldots)$
be the shift on $\Sym$.
We will often use the notation $\underline{s} = (s_1s_2s_3\ldots)$. For all $z
\in \C$ with $E^{\circ n}(z) \in \Cp \mbox{ for all } n \in
\N$ the {\em external address} $S(z) \in \Sym$ is the sequence of the symbols
of the strips containing $z,E(z), E^{\circ 2}(z),\ldots$
\end{definition}

Set $|(n,L)|=|(n,R)|:=|n|$ for $n \in \Z$.

\begin{definition}[Minimal Potential of External Addresses]
\label{minpot} \lineclear
For sequence $\su = s_1s_2s_3\ldots \in \Sym$, define its {\em minimal
  potential} $t_{\su} \in \R^+_0 \cup \{\infty\}$ via
\[
t_{\su}=\inf\left\{t>0\colon
   \lim_{k\geq1}\frac{|s_k|}{F^{\circ(k-1)}(t)}=0\right\}.
\]
\end{definition}
Notice that $t_{\sigma(\su)}=F(t_{\su})$.

\begin{definition}[Exponentially Bounded]
\label{expbound}\lineclear
A sequence $\su \in \Sym$ is {\em exponentially bounded} if there is an $ x>0$
such that $|s_k|\leq F^{\circ(k-1)}(x)$ for all $k\geq1$. 
\end{definition}
This condition is preserved under the shift, but the constant changes:
\[
\su'=\sigma(\su) \quad\Longrightarrow\quad |s_k'|\leq
F^{\circ(k-1)}(F(x)).
\]
An equivalent definition of exponential boundedness (as used in
\cite{Escaping}) is the existence of
$x,A>0$ with $|s_k|\leq AF^{\circ(k-1)}(x)$ for all $k\geq 1$.
It was shown in \cite[Theorem~4.2 (1)]{Escaping} that a sequence $\su$ is
exponentially bounded iff $t_{\su}<\infty$.

\begin{lemma}[External Addresses are Exponentially Bounded]
\label{extadrexpbound} \lineclear
For $E=E_{a,b}$ choose $\delta > 0$ with $|a|+|b|\leq
e^\delta-(\delta+1)$. Then every orbit $(z_k)$ satisfies the bound
\[
\max\{|\Re (z_k)|,|\Im (z_k)|\}\leq|z_k|<F^{\circ(k-1)}(|z_1|+\delta).
\]
In particular, every orbit in $\Cp$ has exponentially bounded external address.
\end{lemma}
\proof
For all $k$ we can estimate
\begin{eqnarray*}
|z_{k+1}|+\delta
&=&
   |ae^{z_k} + be^{-z_k}| + \delta \leq |a| \exp (\Re (z_k)) + |b| \exp
   ( \Re (-z_k))+\delta\\
   &\leq&
   (|a|+|b|)\exp(|z_k|)+\delta \leq (e^\delta
-(\delta+1))\exp(|z_k|)+\delta \\
&=& e^{|z_k|+\delta}-(\delta+1)e^{|z_k|}+\delta \leq
  e^{|z_k|+\delta}-(\delta+1)+\delta =F(|z_k|+\delta) \,\,.
\end{eqnarray*}
Induction yields
\[
|z_k|+\delta< F^{\circ(k-1)}(|z_1|+\delta)
\]
for all $k\ge 1$. If the orbit avoids $\A$, then the external address 
$s_1s_2s_3
\ldots$ is defined and we have (using Lemma~\ref{LemHeightStrips} in the first
inequality):
\[
2\pi|s_k| \leq | \Im (z_k) |+3\pi\leq |z_k|+\delta +3 \pi - \delta <
F^{\circ(k-1)}(|z_1|+\delta) +3\pi - \delta.
\]
\qed

\section{Tails of Dynamic Rays}
\label{SecRayTails}

\Intro{In this section, we show that the set of escaping points of $E$ contains
uncountably many curves starting at $\infty$; in Section~\ref{SecRays}, these
curves will be extended further, and in
Section~\ref{SecClassification}, we will 
show that all escaping points are associated to these curves.}

\begin{definition}[Tail of Ray]
\label{tail} \lineclear
A {\em ray tail} with external address $\su \in \Sym $ is an injective curve
$g_{\su}\colon [\tau,\infty) \rightarrow \C \quad (\tau>0)$ with the
following properties:
\begin{itemize}
\item each point on the curve escapes within $\Cp$
\item each point on the curve has external address $\su$
\item either $\lim_{t\rightarrow \infty} \Re (g_{\su}(t))= +\infty$ or
  $\lim_{t\rightarrow \infty} \Re (g_{\su}(t))= -\infty$ 
\item all $t \geq \tau$ satisfy $E^{\circ k}(g_{\su}(t))=\pm F^{\circ
     (k-1)}(t) +2\pi is_k + O(1)$ as $k\rightarrow \infty$.
\end{itemize}
\end{definition}
The real value $t$ is called the {\em potential} of $z=g_{\su}(t)$.

Given an exponentially bounded external address $\su$, let 
$T_\su\geq T_{a,b}$ be the least value for which 
$4\pi |s_{n+1}| < F^{\circ n}(T_\su)$ for all $n\in\N$ (such a finite $T_\su$
exists since $\su$ is exponentially bounded). Note that 
\(
T_{\sigma(\su)}\leq F(T_{\su})
\).

\begin{lemma}[Minimal Potentials]
\label{LemMinimalPotentials} \lineclear
For every exponentially bounded $\su$ and every $t>t_{\su}$, there is an
$N\in\N$ such that for all $n\geq N$, $F^{\circ n}(t)>T_{\sigma^{n}(\su)}$.
\end{lemma}
\proof
Since $t>t_\su$, there is an $N>0$ such that $F^{\circ n}(t)>4\pi|s_{n+1}|$
for all $n\geq N$; enlarge $N$ if necessary so that also 
$F^{\circ N}(t) > T_{a,b}$. Then 
$F^{\circ n}(t)=
F^{\circ(n-N)}(F^{\circ N}(t))>F^{\circ(n-N)}(T_{\sigma^{N}(\su)})
\geq T_{\sigma^{ n}(\su)}$.
\qed

\begin{proposition}[Existence of Tails of Rays]
\label{extail} \lineclear
For all  $ a,b \in \Cstar$ and for every exponentially bounded sequence $\su
\in \Sym$ there exists a ray tail $g_{\su}\left([T_{\su},
\infty[\right)$ with external address $\su$. 
Furthermore
\[ 
g_{\su}(t)=
\left\{
\begin{array}{ll}
t-\alpha + 2 \pi i s_1 +r_{\su}(t)
& \mbox{(if $s_1\in\Z_R$)} \\
-t+\beta+ 2 \pi i s_1 +r_{\su}(t)
& \mbox{(if $s_1\in\Z_L$)}
\end{array}\right.
\]
with  
\[
|r_{\su}(t)|\leq (C_1+8\pi |s_2|) e^{-t}
\]
with a constant $C_1$ depending only on $a$ and $b$.
Moreover, for $t\geq T_\su$, 
\[
E(g_{\su}(t))=g_\sigsu(F(t)) \,\,.
\]
\end{proposition}
\overview
We construct a family of maps $g_{\su}^n\colon\R^+\to\C$ for $n \in \N$ as
follows:
\begin{equation}
g_{\su}^n(t):=L_{s_1}\circ \ldots \circ L_{s_n}\circ \left( \pm F^{\circ
     n}(t) +2\pi i s_{n+1}\right) \,\,.
\label{EqDefGns}
\end{equation}
The $\pm$ depends on $s_{n+1}$: it is $+$
for $s_{n+1} \in \Z_R$ and $-$ for $s_{n+1} \in \Z_L$ 
(then $\pm F^n(t) + 2\pi i s_{n+1} \in \overline{R_{s_{n+1}}}$ for
large $t$).
We will show that the $g_\su^n$ are defined
for all $t \geq T_\su$ independently of  $n$ and converge uniformly to the
desired function $g_{\su}(t)$.
For this proof we need several lemmas. The first of them gives control
over the inverse branches of $E$.

\begin{lemma}[Control on $L_s$]
\label{LemControlLs} \lineclear
Let $w \in \Cp$ and $z= L_s(w)$ with $|w| > \max \left\{
\sqrt{\left|\frac{2b}{a}\right|}(|a|+|b|)
,\sqrt{\left|\frac{2a}{b}\right|}(|a|+|b|),8  |ab|,1\right\}$ and
$s \in \Z_R$.  Then  
\[
z=\ln w - \alpha +2\pi i p +r^* 
\]
for some $p\in \Z$ with $|r^*|< |2\frac{b}{a}e^{-2z}|<1$
and $|r^*|<8  |ab|\cdot |w|^{-2}<\frac{1}{|w|}<1$.

Similarly, for $s \in \Z_L$ we get 
\[
z=-\ln w +\beta +2 \pi i p +r^* 
\]
for some $p \in \Z$ with $|r^*|<|2\frac{a}{b}e^{-2z}|<1$ and $|r^*|< 8 
|ab|\cdot |w|^{-2}<\frac{1}{|w|}<1$.
\end{lemma}
Note that the branch of $\ln w$ in this lemma is immaterial because of the
ambiguity in $p$.

\proof
Heuristically, $E(z)\approx ae^z$ if $\Re(z)\gg 0$ and $E(z)\approx be^{-z}$ if
$\Re(z)\ll 0$. We discuss the case $s\in \Z_R$; the other case is similar.
We have

\[
w=E(z)=ae^z+be^{-z}=ae^z\left(1+\frac{b}{a}e^{-2z}\right) 
\]
\begin{eqnarray}
\Longrightarrow\quad 
z\in \ln\frac{w}{a} - \ln\left(1+\frac{b}{a}e^{-2z}\right) +2 \pi i
\Z \label{eq1}
\end{eqnarray}

Lemma~\ref{LemSignRealPart} implies $\Re(z)>0$, so we have
\[
    \sqrt{\left|\frac{2b}{a}\right|}(|a|+|b|)    
< |w| \leq |a|e^{\Re(z)}+|b|e^{-\Re(z)} \leq (|a|+|b|)e^{\Re(z)} \,\,,
\]
hence $e^{\Re(z)} > \sqrt{|2b/a|}$ and thus
\[
\Rightarrow \left| \frac{b}{a}
e^{-2z} \right|=\left|\frac{b}{a}\right|e^{-2\Re(z)} < \frac{1}{2}.
\]

Since $|\ln(1+u)| < (2\ln2)|u|$ for $|u| < \frac{1}{2}$, it follows that
\[
\left|\ln\left(1+\frac{b}{a}e^{-2z}\right)\right| 
<2\ln 2 \left|\frac b a e^{-2z}\right|<1
\]
(here we use the principal branch of $\ln$).

Together with (\ref{eq1}) the first claim follows.
For the estimate of $|r^*|$ in terms of $|w|$ it follows from
$2\ln 2 \left|\frac b a e^{-2z}\right|<1$ that
$|be^{-z}|<\frac 1 {2\ln 2}|ae^z|$ and hence by
(\ref{EqStandardEstimate}) we get $ |w|<2|ae^z|$. We thus obtain 
\[
|ab|\frac 2 {|ae^z|^2} < 8 |ab| \cdot |w|^{-2}<|w|^{-1}<1.  
\]
\qed

The next lemma gives us control on $g_\su^1=L_{s_1}\circ(\pm F+2\pi i s_2)$.
\begin{lemma}[Control on $L_s\circ F$]
\label{LemControlLsF} \lineclear
Choose $a,b \in \Cstar$ and $n\in\Z$. Then for $t>0$ such that
\[F(t) > \max\left\{ \sqrt{\left|\frac{2b}{a}\right|}(|a|+|b|),
\sqrt{\left|\frac{2a}{b}\right|}(|a|+|b|),  8  |ab|,1\right\} \]
and $F(t)>4\pi |n|$,

\[
L_{s}(\pm F(t) +2\pi i n)\in t - \alpha + \pi i \Z +r \quad \mbox{ 
with }  \quad |r|<(4+8\pi|n|)e^{-t}<4
\]
 for $s \in \Z_R$ and
\[
L_{s}(\pm F(t) +2\pi i n)\in -t + \beta +  \pi i \Z +r \quad 
\mbox{ with } \quad |r|<(4+8\pi|n|)e^{-t}<4
\]
for $s \in \Z_L$.
\end{lemma}

\proof 
For $s \in \Z_R$, Lemma~\ref{LemControlLs} gives
\begin{eqnarray*}
L_{s}(\pm F(t) +2\pi i n)  &\in& \ln\left(\pm F(t)\left(1\pm\frac{2\pi i
n}{F(t)}\right)\right) -\alpha +2\pi i \Z +r^* \\
&=&t+\ln(1-e^{-t})+\ln\left( 1 \pm \frac{2\pi i 
n}{F(t)}\right)-\alpha +2\pi i \Z  +r^* +(i\pi)
\end{eqnarray*}
with $|r^*|<|F(t)+2\pi i n|^{-1}\leq F(t)^{-1} =e^{-t}\frac{e^t}{e^t-1}<2e^{-t}<1$.
Here, the last term $(+i \pi)$ only occurs in the case
$L_s(-F(t) \ldots)$. 

Since $2\pi |n|/F(t)<1/2$, we get 
\[
\left|\ln\left(1\pm\frac{2\pi i n}{F(t)}\right)\right|
\leq (2\ln 2) \frac{2\pi|n|}{F(t)}\leq 8 \pi |n| e^{-t}<2
\]
and
\[
|\ln(1-e^{-t})| < 2e^{-t}<1.
\]
Thus we get
\[
L_{s}(\pm F(t) +2\pi i n)\in t - \alpha +  \pi i \Z +r
\]
with $|r|<(4+8\pi |n|)e^{-t} <4$. 

The proof for $s \in \Z_L$ is analogous.
\qed

\begin{lemma}[Bound on Real Parts]
\label{LemBoundRealParts} \lineclear
For all $\su \in \Sym$, $n\in \N$ and $t>0$ with $t\geq T_{a,b}$ 
\[
\left|\Re\left(g_{\su}^n(t)\right) \right| > t-(M+2)\,\,.
\]
Moreover, if $|\Re(z)|>\max\left\{\frac 1 2 \ln \left|\frac
{2b}{a}\right|,\frac 1 2 \ln \left|\frac {2a}{b}\right|,\ln\frac 4
{|a|},\ln \frac 4 {|b|}\right\} $, 
then
\[
|E'(z)|>\frac 1 2 K e^{|\Re(z)|}>2 \,\, .
\]

\end{lemma}
\proof We will prove the first part via induction simultaneously for all
external addresses $\su$ and all $t\geq T_{a,b}$. For
$n=0$ clearly
\[
\left| \Re(g_{\su}^0(t))\right|=t>t-(M+2)\,\,.
\]
Assume that $\left|\Re(g_{\su}^{n-1}(t))\right| > t-(M+2)$. Then by
Lemma~\ref{LemControlLs} and inductive hypothesis we get
\begin{eqnarray*}
|\Re (g_{\su}^{n}(t))| &=& 
\left|\Re\left(\pm\ln\left(g_{\sigma(\su)}^{n-1}(F(t))\right) \pm
\{\alpha,\beta\} +2\pi ip +r^*\right)\right|  \\ &\geq&
\left|\Re\left(\ln\left(g_{\sigma(\su)}^{n-1}(F(t)\right)\right)\right|
-M -1
=
\left|\ln\left|g_{\sigma(\su)}^{n-1}(F(t))\right|\right|-M-1 \\
&\geq& 
\left|\ln \left|F(t)-(M+2) \rule{0pt}{12pt}
\right|\right|-M-1=\ln(e^t-M-3)-M-1\\
&>&
t-M-2.
\end{eqnarray*}
The last step needs an elementary calculation based on $F(t)\geq
F(T_{a,b})\geq T_{a,b}+M+4$.

For the second part, we write: 
\[
|E'(z)|=\left\{|a|e^{|\Re(z)|} \cdot \left|1 -
  \frac{b}{a} e^{-2|\Re(z)|}\right| \quad \mbox{ if }
   \Re(z)>0 \atop |b|e^{|\Re(z)|} \cdot  \left|1 -
   \frac{a}{b} e^{-2|\Re(z)|} \right| \quad \mbox{ if }
   \Re(z)<0 \right. \,\,.
\]
By hypothesis, $\left|\frac{b}{a} e^{-2|\Re(z)|}\right| < \frac 1 2$, resp.\
$\left|\frac{a}{b} e^{-2|\Re(z)|}\right| < \frac 1 2$, so we get 
\[
|E'(z)|>\frac 1 2 \min\{|a|,|b|\} e^{|\Re(z)|}
\]
and hence $|E'(z)|>2$.
\qed
We can now finish the construction of dynamic ray tails.

\proofof{\ref{extail}}
We show first that $g_{\su}^n(t)$ converges uniformly in $t$ to a
limit function $g_{\su}([T_\su,\infty[)\to\C$.
For $t \geq T_\su$, we write
\begin{eqnarray*}
 |g_{\su}^{n+1}(t)-g_{\su}^{n}(t)|  &=& 
\lefteqn{\left|L_{s_1}\circ \ldots
\circ L_{s_n}\circ \overbrace{L_{s_{n+1}}\left(\pm F(F^{\circ n}(t))+ 2 \pi i
     s_{n+2}\right)}^{=:\mu_n}\right.}  \\ & & \left.  -L_{s_1}\circ 
\ldots \circ L_{s_n} \circ
\underbrace{\left(\pm F^{\circ n}(t)+2 \pi i s_{n+1}\right)}_{=:\nu_n}\right|
\,\,.
\end{eqnarray*}
Since $\mu_n,\nu_n\in \ovl R_{s_{n+1}}$, we get $|\Im(\mu_n-\nu_n)|\leq 3\pi$ 
by  Lemma~\ref{LemHeightStrips}; by Lemma~\ref{LemControlLsF}, we get
$|\Re(\mu_n-\nu_n)|<4+M$, hence $|\mu_n-\nu_n|<4+3\pi+M$. 

By construction, all $w\in\A$ have $|\Re(w)|\leq 2\sqrt{|ab|}$. However, by
Lem\-ma~\ref{LemBoundRealParts} and (\ref{EqInequTab}),
\begin{eqnarray*}
|\Re\left(L_{s_k}\circ \ldots \circ L_{s_n}(\nu_n) \right)|
&>&
F^{\circ(k-1)}(t)-(M+2) \geq T_{\su}-(M+2) 
\\
&>&
\max\{8|ab|,1\}
\geq
2\sqrt{|ab|}
\end{eqnarray*}
and similarly for $\mu_n$, so the same branch of $L_{s_{k-1}}$ applies to
$L_{s_k}\circ \ldots \circ L_{s_n}(\nu_n)$ and $L_{s_k}\circ \ldots \circ
L_{s_n}(\mu_n)$. 

By Lemma~\ref{LemBoundRealParts},
$L_{s_k}\circ \ldots \circ L_{s_n}(\nu_n)$ and $L_{s_k}\circ \ldots
\circ L_{s_n}(\mu_n)$ are both in the domain
$\left\{z\in\C\colon \Re(z)>T_{\su}-(M+2) \right\}$ or both in
$\left\{z\in\C\colon \Re(z)<-T_{\su}+(M+2) \right\}$ on which
$|E'(z)|>\frac 1 2 Ke^{|\Re(z)|}>2$.
In particular, 
\[
|\Re(\nu_n)|=F^{\circ n}(t)>t+(M+4)
\qquad\mbox{and so}\qquad
|\Re(\mu_n)|>t \,\,.
\]
Therefore, 
\[
|L_{s_n}(\mu_n)-L_{s_n}(\nu_n)| < \frac{|\mu_n-\nu_n|}{\frac1 2 K e^t}
< \frac{8+6\pi+2M}{K}e^{-t}
\,\,.
\]
After repeated application, we get
\begin{eqnarray}
|g_{\su}^{n+1}(t)-g_{\su}^{n}(t)| 
<\frac{8+6\pi+2M}{2^{n-1}K} e^{-t} \,\,.
\label{EqTailContraction}
\end{eqnarray}
Therefore, the $g_\su^n$ converge uniformly to a continuous limit function
$g_\su\colon[T_\su,\infty)\to\C$.

By construction, $E(g_\su^n(t))=g_{\sigsu}^{n-1}(F(t))$, so in the limit we
obtain the desired relation
\[
E(g_\su(t)) = g_{\sigsu}(F(t))
\]
for $t \geq T_\su$. In order to estimate $r$, we get
\begin{eqnarray*}
|g_{\su}(t)-g_{\su}^1(t)|
&\leq&
\sum_{n=1}^{\infty}|g_{\su}^{n+1}(t)-g_{\su}^{n}(t)| \\
&<&
\sum_{n=1}^{\infty} 
\frac{8+6\pi+2M}{2^{n-1}K} e^{-t}
= 2\frac{8+6\pi+2M}{K} e^{-t} \,\,.
\end{eqnarray*}

By Lemma~\ref{LemControlLsF}, we have $g_{\su}(t)=t-\alpha+2\pi i s_1
+r_{\su}(t)$, resp. $g_{\su}(t)=-t+\beta +2\pi i s_1 +r_{\su}(t)$ with
\[
r_{\su}(t) \leq 
\left(2\frac{8+6\pi +2M}{K} +4 + 8\pi |s_2|\right) e^{-t}.
\]

Finally we prove injectivity.
If $g_{\su}(t_1)=g_{\su}(t_2)$ for $t_2 \ge t_1 \geq T_\su$, then we get
\[
 g_{\sigma^n(\su)}(F^{\circ n}(t_1))=g_{\sigma^n(\su)}(F^{\circ n}(t_2))
\]
for all $n \geq 0$ with the bounds 
\begin{eqnarray*}
&&\left|g_{\sigma^n(\su)}(F^{\circ n}(t_j)) - 
\left(\pm F^{\circ n}(t_j) + \{-\alpha
\mbox{ or} +\beta\} -2\pi i s_{n+1}\right)\right|
\\
&=&
|r_{\sigma^n(\su)}(F^{\circ
n}(t_j))|< (C_1+8\pi|s_{n+2}|)e^{-F^{\circ n}(t_j)} \\
&=&
8\pi|s_{n+2}|/F^{\circ(n+1)}(t_j) + o(1)
\,\,.
\end{eqnarray*}
Since $t_1,t_2>t_{\su}$, the right hand side is bounded; this implies that 
$|F^{\circ n}(t_1)-F^{\circ n}(t_2)|$ must be bounded as well as $n\to\infty$.
This implies that $t_1=t_2$.
\qed

\section{Dynamic Rays}
\label{SecRays}

\Intro{
In this section we construct dynamic rays by extending ray tails to as low
potentials as possible. The idea is to use the relation
$E(g_{\su}(t))=g_{\sigma(\su)}(F(t))$ to pull ray tails back by the dynamics:
$g_{\su}$ is a branch of $E^{-1}\circ g_{\sigma(\su)}\circ F$.
}

\begin{theorem}[Existence of Dynamic Rays]
\label{ray} \lineclear
(1) If neither of the two critical orbits escapes then for every exponentially
bounded $\su$ there exists an injective curve 
$g_{\su}:]t_{\su},\infty[ \rightarrow \C $ consisting of escaping points
such that 
\[
E(g_{\su}(t))=g_{\sigma(\su)}(F(t)) \mbox{ for all } t> t_{\su} 
\]
which extends the ray tail with external address $\su$ as constructed in
Proposition~\ref{extail}. In particular, it inherits its asymptotics for large
$t$. 

\noindent
(2) If at least one of the critical orbits escapes, then (1) is still true
for every $\su$, unless $\su$ is such that there is an $n \geq 1$ and a $t_0 >
F^{\circ n}(t_{\su})$ with $ g_{\sigma^n(\su)}(t_0) \in\{v_1,v_2\}$. For
those exceptional $\su$, there is an injective curve
$g_{\su}:]t^*_{\su},\infty[\rightarrow \C$ with the same properties as
before, where $t^*_{\su}$ is the largest potential which has an $n
\geq 1$ such that $g_{\sigma^n (\su)}(F^{\circ n}(t^*_{\su}))\in\{v_1,v_2\}$.
\end{theorem}

The curve $g_{\su}:]t_{\su},\infty[ \rightarrow \C $ is called the {\em 
dynamic ray at external address $\su$}.

\remark
Notice that it is no longer required that all points on the ray
share the external address $\su$. Since the partition $\{R_s\}$ of $\C$ is
unnatural from the dynamical point of view, there is no reason why the rays
should respect it. Only the points with large potential (those which are on ray
tails) have external address $\su$. Note that the definition of $t_\su$ is the
same as for the case of exponential functions $\lambda\exp$ in \cite{Escaping}:
since $E_{a,b}(z)\approx ae^z$ resp.\ $E_{a,b}(z)\approx be^{-z}$ in far right
resp.\ left half planes, this coincidence of the values of $t_\su$ reflects the
fact that the value of $t_\su$ for exponential functions does not
depend on the complex parameter $\lambda$.

\proof 
We want to show that $g_{\su}(t)$ exists for $t>t_{\su}$. Choose $\eps>0$. 
By Lemma~\ref{LemMinimalPotentials}, there is an $N\in\N$ such that
$F^{\circ N}(t_\su+\eps)>T_{\sigma^{N}(\su)}$.

By Proposition~\ref{extail}, the ray tail
$g_{\sigma ^N(\su)}(t)$ exists for $t\geq T_{\sigma^{N}(\su)}$. 
Since $E^{\circ N}$ maps the ray tail $g_\su$ to a tail of
$g_{\sigma^N{\su}}$, it follows that there is a branch of $E^{-N}$ which sends
(a subset of) the ray tail $g_{\sigma^{N}}\colon [T_{\sigma^{N}(\su)},\infty)\to\C$ to a curve $g_\su\colon [t_\su+\eps,\infty)\to\C$
extending the ray tail $g_\su\colon[T_\su,\infty]$. In other words, the ray
tail $g_\su$ can be extended to potentials $t_\su+\eps$. Since $\eps$ was
arbitrary, we have shown the existence of the dynamic ray
$g_{\su}\colon(t_\su,\infty)\to\C$. Injectivity follows from
Proposition~\ref{extail}. This proves the second claim.

Note that the pull-back of a ray tail is possible if and only if it contains no
critical value. Therefore, this construction can be carried out for all ray
tails except those mentioned in the exceptions of the third claim. This proves
the theorem.
\qed

\remark
Note that the pull-back in the proof of the previous theorem need not respect
the partition we had initially used for constructing the ray tails: the branch
of every $E^{-1}$ is determined using the ray tails and is then continued
analytically. It may well happen that a dynamic ray crosses the partition
boundary, but only at potentials below $T_\su$.

Next we investigate under which conditions an escaping point is on a ray.

\begin{theorem}[Fast Escaping Points are on Ray]
\label{fast} \lineclear
Let $(z_k)$ be an escaping orbit and let $\su$ be such that $z_k \in
\overline{R_{s_k}}$ for all $k$; suppose also that there exists a $t'
> t_{\su}$ with $|\Re(z_k)|\geq F^{\circ (k-1)}(t')$ for infinitely many $k$.
If in addition 
\begin{equation}
|\Re(z_k)|>2\sqrt{|ab|}+1
\quad\mbox{for all $k$, as well as}
\quad
t'\geq T_{\su}
\label{Eq z_k in regular domain}
\end{equation}
then there is a $t\geq t'$ such that $z_1=g_{\su}(t)$. 

If (\ref{Eq z_k in regular domain}) is not satisfied for all $k$, then there is
an $N\in\N$ such that $z_{N+1}=g_{\sigma^{N}(\su)}(F^{\circ N}(t))$ for
some $t\geq t'$; moreover, at least if the two critical orbits do not escape,
there is an external address $\su'$ which differs from $\su$ only in
finitely many entries such that $z_1=g_{\su'}(t)$ for some $t\geq
t'>t_{\su'}=t_{\su}$. 
\end{theorem}
\proof
For $k \in \N$ pick $t_k >0$ such that $F^{\circ (k-1)}(t_k)=|\Re(z_k)|$. 
By assumption, $t_k \geq t'$  for infinitely many $k$, and by
Lemma~\ref{extadrexpbound}, the sequence $(t_k)$ is bounded above. 
Moreover, by Lemma~\ref{LemHeightStrips}, we have $|\Im(z_k)| \leq
2\pi(|s_k|+3/2) \ll F^{\circ(k-1)}(t')$ for all sufficiently large $k$. 
Therefore, $|\Im(z_k)| \ll F^{\circ
(k-1)}(t')\leq |\Re(z_k)|$ for infinitely many $k$.
From now on we will only look at such $k$.

Suppose first that $|\Re(z_k)|>2\sqrt{|ab|}+1$ for all $k$. Then all
$z_k\notin\A$, hence $z_k\in R_{s_k}$, and thus
\[
z_l=L_{s_l}\circ\dots\circ L_{s_{k-1}}(z_k) 
\]
for all $l\in\{k-1,\dots,2,1\}$. 

Consider the points $w_k:=\pm F^{\circ(k-1)}(t_k)+2\pi i s_k$, where the
signs are chosen such that for all $k$, $z_k$ and $w_k$ are in the same (right
or left) half plane; then, using Lemma~\ref{LemHeightStrips} again,
\begin{equation}
\left|z_k - w_k\right|
\leq 3 \pi.
\label{Eq z_k Diff}
\end{equation}
Therefore, given $\eps\in]0,1[$, 
the derivative bound in Lemma~\ref{LemBoundRealParts} implies
\(
\left|z_{k-1}-L_{s_{k-1}}(w_k)\right|< \eps 
\)
for sufficiently large $k$.
But this implies that the same branch of $L_{s_{k-2}}$ applies to both
points, and repeated application of this argument shows that if $k$ is
sufficiently large, then for all $l\in\{k-1,\dots,1\}$ we have
\[
\left| z_l-L_{s_l}\circ\dots\circ L_{s_{k-1}}
\left(w_k\right)\right| <\eps
\,\,.
\]
If in addition $t_k\geq T_{\su}$, then $g_{\su}^{k-1}(t_k)
=L_{s_1}\circ\dots\circ L_{s_{k-1}}
\left(w_k\right)$, and we have $|z_1-g_{\su}^{k-1}(t_k)|<\eps$.
Now suppose $t'\geq T_{\su}$. 
Let $t$ be a limit point of the sequence $(t_k)$ (restricted to such $k$ as
mentioned above).
  Obviously $t \geq t' \geq T_{\su}$. By uniform convergence of
  $g_{\su}^{k-1}$ to $g_{\su}$ for potentials at least $T_{\su}$, we get 
$|g_{\su}^{k-1}(t_k)-g_{\su}(t_k)|< \eps$ (possibly by enlarging $k$).
Finally, for $t_k$ close enough to $t$ ($t$ is a limit point of $(t_k)$)
$|g_{\su}(t_k)-g_{\su}(t)|< \eps$.
Combining this, it follows
\[
|z_1-g_{\su}(t)|\leq |z_1 - g_{\su}^{k-1}(t_k)|+|g_{\su}^{k-1}(t_k) -
g_{\su}(t_k)| + |g_{\su}(t_k) - g_{\su}(t)| < 3\eps
\]
for certain sufficiently large $k$. Hence $g_{\su}(t)=z_1$ because
$\eps>0$ was arbitrary.

If the condition $|\Re(z_k)|>2\sqrt{|ab|}+1$ does not hold for all $k$, then
there is an $N\in\N$ such that it holds for all $k\geq N$; similarly, by
Lemma~\ref{LemMinimalPotentials}, if $t'\geq T_{\su}$ is not satisfied, then
for sufficiently large $N$, we have $F^{\circ N}(t')\geq T_{\sigma^{N}(\su)}$.  Therefore, there is an $N\in\N$ with $z_{N+1}=g_{\sigma^{N}(\su)}(F^{\circ N}(t))$ for some $t\geq t'$. Pulling back $N$ times along
the orbit from $z_1$ to $z_{N+1}$, it follows that $z_1=g_{\su'}(t)$ for some
external address $\su'$ which can differ from
$\su$ only in the first $N$ entries; thus $t_{\su'}=t_{\su}$. However, if
the pull-back runs through a dynamic ray which contains an escaping critical
value, then this pull-back is impossible --- and only then.
\qed

\begin{proposition}[Controlled Escape for Points on Rays]
\label{control} \lineclear
For every exponentially bounded external address $\su$ and for every
$t>t_{\su}$, the orbit of $g_{\su}(t)$ satisfies the asymptotic bound
\[
E^{\circ k}(g_{\su}(t)) = 
\left\{\begin{array}{ll}
F^{\circ k}(t) - \alpha  +2\pi i s_{k+1} +o(1)
& \mbox{(if $s_{k+1}\in\Z_R$)} \\
-F^{\circ k}(t)  + \beta +2\pi i s_{k+1} +o(1)
& \mbox{(if $s_{k+1}\in\Z_L$)}
\end{array}\right.
\]
as $k\to \infty$.
In particular, for every real $p>0$ it satisfies 
\[
\frac{|\Im (E^{\circ k}(g_{\su}(t)))|^p}{\Re(E^{\circ
 k}(g_{\su}(t)))}\longrightarrow 0 \,\,.
\]
\end{proposition}
\proof
By Proposition~\ref{extail}  we have good error bounds for the dynamic rays
$g_{\sigma^k(\su)}$  for potentials greater than $T_{\sigma^{k}(\su)}$.
Since $t>t_\su$, there exists a $k_0$ such that $F^{\circ
k}(t)>T_{\sigma^{k}(\su)}$ for all $k \geq k_0$
(Lemma~\ref{LemMinimalPotentials}). Then
\begin{eqnarray*}
E^{\circ k}(g_{\su}(t))=g_{\sigma^k(\su)}(F^{\circ k}(t))&=&F^{\circ
   k}(t)-\alpha+2 \pi i s_{k+1}+r_{\sigma^k(\su)}(F^{\circ k}(t))\\
\mbox{ resp.  } &=&-F^{\circ
   k}(t)+\beta+2 \pi i s_{k+1}+r_{\sigma^k(\su)}(F^{\circ k}(t))
\end{eqnarray*}
with
\[
|r_{\sigma^k(\su)}(F^{\circ k}(t))|< (C_1+8\pi|s_{k+2}|)e^{-F^{\circ k}(t)},
\]
where $C_1$ only depends on $a$ and $b$.
This tends to $0$ as $k \rightarrow \infty$ (extremely fast).
Along the orbit of $g_{\su}(t)$, the absolute values of real parts thus grow
like $F^{\circ k}(t)$, while the imaginary parts are bounded in absolute value
by the asymptotically much smaller quantity $AF^{\circ k}(t_{\su}+\eps)$ for
any $\eps\in]0,t-t_{\su}[$. In particular we have
\[
\ln\left(\frac{|\Im (E^{\circ k}(g_{\su}(t)))|^p}{|\Re(E^{\circ
k}(g_{\su}(t)))|}\right) < pF^{\circ (k-1)}(t_{\su}+\eps)-F^{\circ (k-1)}(t)
+O(1) \rightarrow - \infty. 
\]
This proves the last claim.
\qed

\section{Eventually Horizontal Escape}
\label{SecHorizontalEscape}

\Intro{The main result in this section shows essentially that, for any given
external address, there is at most one point with this external address which
is not on a dynamic ray.}

Define
\[
R_{h}:=\max \left\{ 
\ln\left(\frac{2h+8\pi}{|a|\pi}\right),
\ln\left(\frac{2h+8\pi}{|b|\pi}\right),
\frac{1}{2}\ln\left|\frac{2b}{a}\right|,\frac{1}{2}\ln\left|\frac{2a}{b}\right|,
|a|+|b|
\right\}
\]
(depending on the parameters $a,b\in\Cstar$).

\begin{lemma}[Exponential Separation of Orbits]
\label{LemExpoSeparation} \lineclear
Let $R > R_{h}$ and let $(z_k)$ and $(w_k)$  be two escaping
orbits with  $|\Re(z_k)|>R$, $|\Re(w_k)|>R$, 
$\Re(z_k)/\Re(w_k)>0$ and $|\Im(z_k-w_k)|<h$ for some
$h>0$ and all $k$. Define
$d_k:=\Re(z_k-w_k)$ for all $k$. If $|\Re(z_1)|- |\Re(w_1)| \geq 3$ and none
of the critical values escapes, then the following holds
\begin{enumerate}
\item $|d_{k+1}| \geq \exp(|d_k|)$ and $|\Re(z_k)|=|\Re(w_k)|+|d_k|$ 
for all $k \geq 1$;
\item 
if $\su\in\Sym$ is such that $z_k\in R_{s_k}$ for all $k$, then
$z_1 = g_{\su'}(t')$ for some $t'> t_{\su}$ and some $\su'$ which differs
from $\su$ in only finitely many entries;
if $|d_1|$ is sufficiently large, then $\su'=\su$;
\item if $w_1 = g_{\su''}(t'')$ for some external address $\su''$ and some
$t''>t_{\su''}$, then $t'>t''$.
\end{enumerate}
If one or both the critical values escape then the first and third statements
are still true, and there is an $N\in\N$ such that $z_{N+1}=
g_{\sigma^{N}(\su)}(F^{\circ N}(t))$ for some $t> t_{\su}$.
\end{lemma}

\proof
1. Let $w_k=t_k +i u_k$ with real $t_k,u_k$. We write
\[ 
|w_k|=|E(w_{k-1})|=\left\{|a|e^{t_{k-1}} \cdot \left|1 +
  \frac{b}{a} e^{-2 w_{k-1}}\right| \quad \mbox{ if }
   t_{k-1}>0 \atop |b|e^{-t_{k-1}} \cdot  \left|1 +
   \frac{a}{b} e^{2 w_{k-1}} \right| \quad \mbox{ if }
   t_{k-1}<0 \right..
\]
By assumption, $\Re(z_k)$ and $\Re(w_k)$ always have the same sign. 
By definition of $d_k$, the property
$|\Re(z_k)|=|\Re(w_k)|+|d_k|$  is equivalent to the fact that $d_k$ has the
same sign as $\Re(z_k)$ and
$\Re(w_k)$, i.e.\ $\frac{t_k}{d_k} >0$.
Since $|\Re(w_{k-1})|>R$, we have 
$e^{-2|\Re(w_{k-1})|} <
\min
\left\{
\left|
\frac{a}{2b} \right|, \left| \frac{b}{2a} \right| \right\}$ and thus
\begin{eqnarray}
\frac{3}{2}|a| e^{|t_{k-1}|}\geq |w_k| \geq \frac{1}{2}|a|
e^{|t_{k-1}|}
\quad
\mbox{if } t_{k-1} >0 \,\,,   
\label{EqInequ_W_k_a}
\end{eqnarray}
and
\begin{eqnarray}
\frac{3}{2}|b| e^{|t_{k-1}|}\geq |w_k| \geq \frac{1}{2}|b|
e^{|t_{k-1}|}
\quad
\mbox{if } t_{k-1} <0 \,\,. 
\label{EqInequ_W_k_b}
\end{eqnarray}
Since $\Re(w_k)$ and $\Re(z_k)$ always have the same signs, there exists
a  $K'\in \{|a|,|b|\}$ (depending on $s_k$) such that 
\begin{equation}
\frac{3K'}{2}e^{|\Re(z_{k-1})|} 
\geq|z_k|\geq
\frac{K'}{2}e^{|\Re(z_{k-1})|}
\label{EqInequ_z_k}
\end{equation}
and 
\begin{equation}
\frac{3K'}{2}e^{|\Re(w_{k-1})|} 
\geq|w_k|\geq \frac{K'}{2}e^{|\Re(w_{k-1})|} \,\,.
\label{EqInequ_w_k}
\end{equation}

Now let $|\Re(z_{k-1})|-|t_{k-1}|\geq 3$ be the inductive hypothesis; then 
$|d_{k-1}|\geq 3$ and $\frac{t_{k-1}}{d_{k-1}}>0$.
Since $|\Im(z_k-w_k)|<h$, the Pythagorean Theorem implies
\begin{eqnarray}
(\Re(z_k))^2&=&|z_k|^2-(\Im(z_k))^2 \geq |z_k|^2 -(|w_k|+h )^2 \nonumber\\
&\geq& \left(\frac{K'}{2} \exp|t_{k-1}+d_{k-1}|\right)^2-\left(\frac{3}{2}K'
\exp|t_{k-1}|+h\right)^2 \,\,.
\label{EqPythagorean}
\end{eqnarray}
This difference is positive: from
$|t_{k-1}| > R_h>\
\max \left\{
\ln\left(\frac{2h}{|a|\pi}\right),\ln\left(\frac{2h}{|b|\pi}\right)
\right\}$, it follows $\frac{2h}{K'e^{|t_{k-1}|}}<\pi$; since 
 $|d_{k-1}|\geq3$, we get
$e^{|d_{k-1}|}\geq e^3 > 3+\pi>3+
\frac{2h}{K'e^{|t_{k-1}|}}$ and therefore 
\[
\frac{K'}{2}e^{|t_{k-1}|}e^{|d_{k-1}|}\geq \frac{3K'}{2} e^{|t_{k-1}|}+h.
\]
Hence we can extract the root in (\ref{EqPythagorean}) and obtain
\begin{eqnarray}
|t_k+d_k|\geq \frac{K'}{2}
  \exp|t_{k-1}+d_{k-1}|\sqrt{1-\left(\frac{3\exp|t_{k-1}|+2h/K'}
{\exp|t_{k-1}+d_{k-1}|}\right)^2}. \label{EqPythagorean2}
\end{eqnarray}
Since $|d_{k-1}|\geq 3$ and $\frac{t_{k-1}}{d_{k-1}}>0$, we get as above
\begin{eqnarray}
0&<& \frac{3\exp|t_{k-1}|+2h/K'}{\exp|t_{k-1}+d_{k-1}|}
=\frac{3 +\left(2h/K'\right) \exp(-|t_{k-1}|)}{\exp|d_{k-1}|} \nonumber \\
&<& \frac{3+\pi}{\exp|d_{k-1}|}\leq \frac{3+\pi}{e^3} < \frac{1}{2} \,\,.
\label{EqRadicand}
\end{eqnarray}
The radicand in (\ref{EqPythagorean2}) is thus in $\left(\frac{1}{2}, 1 \right)$.
Since $|t_k|=|\Re(w_k)| \leq |w_k|$, we get
\begin{eqnarray*}
|t_k+d_k|-|t_k| &\geq&  \frac{K'}{4} \exp|t_{k-1}+d_{k-1}|
  -\frac{3}{2}K' \exp|t_{k-1}| \\
&=& \frac{K'}{4} \exp|t_{k-1}| \left(\exp|d_{k-1}|-6\strut\right) > 0.
\end{eqnarray*}
Thus $|\Re(z_k)|>|\Re(w_k)|$ and 
$\frac {t_k}{d_k}>0$. Since $|t_{k-1}| \geq  R$, we get $ \frac{K'}{4}
\exp|t_{k-1}| > \frac{K'}{4}\exp(\ln(8/K'))=2$. Thus
\begin{eqnarray*}
|d_k|&=&|t_k+d_k|-|t_k| > 2 (\exp|d_{k-1}|-6)\\
&=& \exp|d_{k-1}|+(\exp|d_{k-1}|-12) \geq \exp|d_{k-1}| \,\,.
\end{eqnarray*}
This shows the first claim.

2. The second claim will be done in four steps: 
\begin{enumerate}
\def\theenumi{(\alph{enumi})}
\item
$|t_k+d_k|\geq F^{\circ(k-N)}(t_{N}+d_{N}+\gamma')-\gamma'$ for some
$\gamma'\in \R$ and $N\in\N$ and for all $k\geq N$ (using the estimates from
the first step)
\label{StepBound_t_k_d_k}
\item
$|t_k|+\gamma\leq F^{\circ(k-N)}(t_{N}+\gamma)$ for some $\gamma\in\R$ and
all $k>N$ (from the maximal growth rate along orbits)
\label{StepBound_t_k}
\item
this will imply an upper bound on $|w_k|$, hence on $|\Im(w_k)|$ and on
$|s_k|$; we will get
$F^{\circ (N -1)} (t_{\su}) \leq |t_{N}| + \gamma$
\label{StepBound_w_k}
\item
the conclusion will then follow.
\label{StepConclusion}
\end{enumerate}

\ref{StepBound_t_k_d_k}
Choose $\gamma' < \ln \left(\frac{K}{4}\right)$ (possibly negative). By
(\ref{EqPythagorean2}), (\ref{EqRadicand}) and since $K'\geq K$, we have
\begin{eqnarray*}
|t_k+d_k|&=&|\Re(z_k)| \geq \frac{K}{4} e^{|t_{k-1}+d_{k-1}|}=
e^{|t_{k-1}+d_{k-1}|+\ln \frac{K}{4}}\\
&\geq&  e^{|t_{k-1}+d_{k-1}|+\gamma'}-\gamma' -1 = 
F(|t_{k-1}+d_{k-1}|+\gamma')-
\gamma'
\end{eqnarray*}
provided $|t_{k-1}+d_{k-1}|$ is sufficiently large. There is thus an $N>0$ such
that inductively for all $k \geq N$
\begin{eqnarray}
|t_k+d_k|+\gamma' \geq F^{\circ (k-N)}(|t_{N}+d_{N}|+\gamma') 
\,\,.
 \label{EqBound_t_k_d_k} 
\end{eqnarray}

\ref{StepBound_t_k}
Choose $\gamma > \ln \left(\frac{3K_{\max}}{2}\right)$. We get
\begin{eqnarray*}
|t_k| &=&  |\Re(w_k)| \leq |w_k| \leq \frac{3K_{\max}}{2} e^{|t_{k-1}|}=
F\left(|t_{k-1}|+\ln
\frac{3K_{\max}}{2}\right)+1 \\ &\leq& F(|t_{k-1}| + \gamma) -\gamma,
\end{eqnarray*}
again provided $t_{k-1}$ is sufficiently large. Possibly by enlarging $N$, we
have for all $k\geq N$
\[
|t_k|+\gamma \leq F^{\circ (k-N)}(|t_{N}|+\gamma) \,\,.
\]

\ref{StepBound_w_k}
For every $\eps >0$ the definition of $t_{\su}$ shows that there are infinitely
many $k\geq N$ with
\begin{eqnarray*}
F^{\circ (k-N)}(F^{\circ (N-1)}(t_{\su} -\eps)) &=& F^{\circ 
(k-1)} (t_{\su}-\eps)
< 2\pi |s_k| \leq | \Im (w_k)| +3 \pi
\\ &\leq& |w_k| +3\pi
\leq
\frac{3K_{\max}}{2} e^{|t_{k-1}|} +3\pi  < e^{|t_{k-1}|+\gamma}+ 3\pi  \\ &=&
F(|t_{k-1}|+\gamma) +1 + 3\pi  \\
  &\leq& F^{\circ(k-N)}(|t_{N}|+\gamma)+1+ 3\pi \,\,.
\end{eqnarray*}
Since this holds for arbitrarily large $k$, we get
\begin{eqnarray*}
F^{\circ (N -1)} (t_{\su}-\eps) \leq |t_{N}| + \gamma \,\,.
\end{eqnarray*}
Since $\eps>0$ was arbitrary, it follows
\[
t_{\sigma^{N-1}(\su)}=F^{\circ(N-1)}(t_\su)\leq |t_{N}|+\gamma \,\,.
\]

\ref{StepConclusion}
Enlarge $N$ if necessary so that
$|t_{N}+d_{N}|+\gamma'>|t_{N}|+\gamma+2$ and $F^{\circ(N-1)}(t_{\su}+2)>
T_{\sigma^{N-1}(\su)}+1$ (such $N$ exists by
Lemma~\ref{LemMinimalPotentials}). Then
\begin{eqnarray*}
|\Re(z_k)|&=&|t_k+d_k|
\geq
F^{\circ(k-N)}\left(|t_{N}+d_{N}|+\gamma'\right)-\gamma'
\\
&>&
F^{\circ(k-N)}\left(|t_{N}|+\gamma+2\right)-\gamma'
\geq F^{\circ(k-N)}\left( t_{\sigma^{N-1}(\su)} )
+2\right)-\gamma' \\
&\gg& F^{\circ(k-N)}\left( t_{\sigma^{N-1}(\su)}+1 \right)
\end{eqnarray*}
for infinitely many $k>N$.
Now Theorem~\ref{fast} implies that there is a $k\geq N$ and a
$t^*>t_{\sigma^{k-1}(\su)}$ such that
\[
z_{k+1}=g_{\sigma^{k}(\su)}(t^*) \,\,.
\]
If no critical value escapes, then we can pull back $k$ times, and it
follows that there is an $\su'$ which differs from $\su$ only in finitely
many entries, and a $t'>t_\su=t_{\su'}$ with
\[
z_1=g_{\su'}(t') \,\,.
\]

If $|d_1|$ is large enough so that $|t_k+d_k|\geq
\max\{F^{\circ(k-1)}(T_{\su}),2\sqrt{|ab|}+1\}$ for all $k$, then the
conditions of Theorem~\ref{fast} are satisfied immediately so that
$z_1=g_{\su}(t)$ for some $t \geq T_{\su}$.

3.
If $w_1 = g_{\su''}(t'')$ for $t''> t_{\su''}$, then $t>t''$ by
Proposition~\ref{control}. 
\qed

\section{Classification of Escaping Points}
\label{SecClassification}

\Intro{
In this section we show that all escaping points are organized in the
form of dynamic rays which are associated to exponentially bounded external
addresses, and we complete the classification of escaping points.
}

\begin{definition}[Limit Set, Landing Point, Uniform Escape]
\label{limits} \lineclear
The {\em limit set of the ray $g_{\su}$} is defined as the set of all possible
limit points of $g_{\su}(t_k) $ for $t_k \searrow t_{\su}$. We  say that 
{\em the ray $g_{\su}$ lands at a point $w$} if $\lim _{t'\searrow
   t_{\su}} g_{\su}(t')$ exists and is equal to $w$ (the limit set consists of
 only one point). If $g_{\su}$ lands at an escaping point $w=g_{\su}(t_{\su})$,
we say that { \em ray and landing point escape uniformly} if for every $R \in
\R$ there exists an $N
\geq 0$  such that for every $n \geq N$, we have $|\Re( g_{\su}^n([t_{\su},
\infty[))|> R$.
\end{definition}

For $R>0$ let $Y_R$ be defined as 
\[
Y_R:=\left\{ z \in \C \colon |\Re(z)|<R \right\} \,\,.
\]

\begin{lemma}[Escaping Set Connected]
\label{LemEscapingConnected} \lineclear
Let $R \geq R_{3\pi}$ with 
\[    
\frac{K}{2} e^R \geq \frac
{2+3\pi}{1-3/e^2}  
\]
and let $(z_k)$ be an escaping orbit which is completely contained in $\C
\sm Y_{R+2}$. Then there exists a closed connected set $C \subset
\Cbar \sm Y_R$ with $\{z_1,\infty\}\subset C$ such that the orbit of every
$z\in C$ is completely contained in $\C \sm Y_R$ and escapes
such that
\[
|\Re(E^{\circ (k-1)}(z))| 
\geq |\Re(z_k)|-2.
\]
Let $\su$ be the external address of $z_1$. All points in $C\sm\{\infty\}$ have
external address $\su$ and 
\begin{itemize}
\item
either $C\sm\{\infty\}=g_\su([t,\infty[)$ for some $t>t_\su$,
\item
or the ray $g_\su$ lands at $z_1$, and
$C=\{z_1,\infty\}\dot\cup g_\su(]t_\su,\infty[)$.
\end{itemize}
In particular, all points in $C\sm\{\infty\}$, except possibly $z_1$, lie on the
 ray $g_{\su}$.
\end{lemma}
\remark
The external address of $z_1$ is defined uniquely because the orbit of $z_1$
never enters the vertical strip $Y_R\supset\A$, so it can never hit the boundary
of our partition.

\proof 
By (\ref{EqInequ_z_k}), for every $k$ there is a
$K'\in\{|a|,|b|\}$ (depending on $k$) such that
\begin{eqnarray}
\frac 3 2 K'\exp|\Re(z_k)|\geq |E(z_k)|= |z_{k+1}|
&\geq& 
\frac{K'}{2}\exp|\Re(z_k)| \geq
\frac{K'}{2}e^R  \nonumber\\
&\geq& \frac{2+3\pi}{1-3/e^2}> 0
\label{EqEscapingC}
\end{eqnarray}
and thus 
\[
|z_{k+1}|-2-3\pi \geq \frac{3}{e^2}|z_{k+1}|.
\]
For $k\geq1$ define
\[
 S_k:= \left\{ \rule{0pt}{24pt} z \in R_{s_k} \quad \colon \left\{ {\Re(z) > \Re(z_k) -2, 
\mbox { if }
       \Re(z_k)>0 \atop \Re(z) < \Re(z_k) +2, \mbox { if } \Re(z_k)<0 }
   \right. \right\}.
\]
Our first claim is $E(S_k) \supset S_{k+1}$ for all $k$. 
Every $z \in R_{s_k} $ with $|\Re(z)|\leq |\Re(z_k)|-2$ satisfies as in
(\ref{EqEscapingC})
\[
|E(z)|\leq
\frac{3}{2}K'\exp|\Re(z)| \leq
\frac{3}{2}K'\exp(|\Re(z_k)|-2)\leq \frac{3}{e^2}|z_{k+1}|.
\]
Since $E: \overline{R_{s_k}} \rightarrow \C$ is surjective, it follows that
$E(S_k)$ contains every point $w \in S_{k+1}$ with 
\[
|w| > \frac{3}{e^2}|z_{k+1}|.
\]
But $w\in S_{k+1}$ implies
\[
|\Re(w)|\geq |\Re(z_{k+1})|-2 \quad \mbox{ and } \quad
|\Im(w)-\Im(z_{k+1})|<3\pi
\]
and thus
\[
|w| \geq |z_{k+1}|-2-3\pi\geq \frac{3}{e^2}|z_{k+1}|.
\]
Hence we get $w \in E(S_k)$ and the first claim is proved. 

Since $\A \cap S_{k+1} = \emptyset$ we obtain a connected set $C'_k
\subset S_k$ such that $E\colon C'_k \rightarrow S_{k+1}$ is a conformal
isomorphism. For $k \geq 1$ consider the sets
\[
C_k:=\left\{z \in S_1 \quad \colon E^{\circ i}(z) \in 
\ovl S_{i+1}, \mbox{ for } i=0,1,\ldots,k-1\right\}\cup\{\infty\} \,\,.
\]
The sets $C_k$ are non-empty (because $\{z_1,\infty\}\subset C_k$), compact and
nested: $C_{k+1}\subset C_k$. 
We just proved that $E^{\circ(k-1)}: C_k \rightarrow
S_k\cup\{\infty\}$ is a homeomorphism for all $k$. Therefore, all $C_k$ are
connected. The nested intersection of non-empty connected compact sets is
non-empty connected and compact; therefore   
\[
C:= \bigcap_{k\geq1} C_k
\]
is a closed connected and compact subset of $\Cbar$ with $\{z_1,\infty\}\subset
C$. Set $C':=C\sm\{\infty\}$. For $z\in C'$ we have $E^{\circ (k-1)}(z) \in
\overline{S_k}$ for all $k$ and thus $|\Re(E^{\circ (k-1)}(z)| > |\Re(z_k)|-2
\rightarrow \infty$. Hence $C'$ only consists of escaping points with orbits in
$\C \sm Y_R$ which have the same external address as $z_1$.

It remains to show that there is at most one point in $C'$ which is
not on the dynamic ray $g_\su$. By Lemma~\ref{LemExpoSeparation}, there is a
$\xi>0$ such that every $z\in C'$ with $|\Re(z)|>\xi$ is on $g_\su$. In
particular, $C'$ contains an unbounded connected part of the tail of $g_\su$.

Suppose that there are two points $z_1,w_1 \in C'$ with orbits $(z_k)$ and
$(w_k)$ such that $|\Re(z_k-w_k)|<3$ for all $k$. By
Lemma~\ref{LemHeightStrips}, we also have $|\Im(z_k-w_k)|<3\pi$, hence
$|z_k-w_k|<3+3\pi$. The derivative $E'$ is bounded below along the orbits of
$z_1$ resp.\ $w_1$ by $|ae^R-be^{-R}|\geq \frac{K}{2}e^R>4$. By pulling back
choosing the branch $E^{-1}: S_{k+1} \rightarrow C_k$, we get
\[
|z_{k-1}-w_{k-1}|<(3\pi+3)/4
\]
and thus inductively $|z_{k-j}-w_{k-j}|<(3\pi+3)/4^j$. Hence 
\(
 |z_1-w_1|<(3\pi +3)/4^{k-1}
\)
for all $k$ and thus $z_1=w_1$. Therefore, if $z_1\neq w_1$, then by
Lemma~\ref{LemExpoSeparation} at least one of these points (say $z_1$)
satisfies $z_1=g_{\su'}(t)$ for an external address $\su'$ which differs from
$\su$ in only finitely many positions, and $t>t_{\su'}=t_\su$; say
$\sigma^{N}(\su')=\sigma^{N}(\su)$. But since $E^{\circ
N}(C')\cup\{\infty\}$ is connected, a single branch of $E^{-N}$ maps $E^{\circ
N}(C')$ to $C'$, avoiding $\A$ in the pull-back, so it follows even that
$\su'=\su$.

Therefore, every point in $C'$ with at most one exception is on $g_\su$.
Since $g_\su\colon]t_\su,\infty[\to\C$ is continuous and $C'$ is closed, the
set $\{t\in\,]t_{\su},\infty[\,\colon g_\su(t)\in C'\}$ is closed in
$]t_\su,\infty[$. We use this to show that if $g_\su(t)\in C'$, then
$g_{\su}(t')\in C'$ for all $t'>t$: otherwise, there would be $t_2>t_1>t_\su$
with $g_\su(t_1),g_{\su}(t_2)\in C'$, but $g_\su(]t_1,t_2[)\cap C'=\emptyset$.
But then at least for large $N$, there would be a large gap between
$g_{\sigma^N(\su)}(F^{\circ N}(t_2))$ and $g_{\sigma^N(\su)}(F^{\circ N}(t_1))$
which could not be filled by points in $g_{\sigma^N(\su)}(]t_{\su},t_1[)$ or
$g_{\sigma^N(\su)}(]t_2,\infty[)$ (Lemma~\ref{LemExpoSeparation}~(3)) or by
the single exceptional point, so $E^{\circ N}(C')$ could not be connected.

Therefore, except for possibly a single point, $C'$ equals either
$g_\su([t,\infty[)$ for some $t> t_\su$, or $g_\su(]t_\su,\infty[)$. But
since $C'$ is closed and connected, an extra point can (and must) occur
only in the second case, and this is what we claimed.
\qed

For an external address $\su$ let $S(\su)$ be the space of external
addresses $\su'$ which differ from $\su$ at only finitely many entries.
Clearly, all $\su'\in S(\su)$ have $t_{\su'}=t_\su$.

\begin{lemma}[Limit Set does Not Intersect Ray]
\label{LemLimitSetNotRay} \lineclear
Let $g_{\su}$ be a ray with the property that the entire orbits of all of its
points avoid $Y_R$ for an $R > 2\sqrt{|ab|}+M+2$ and $R>R_h$ for $h=2\pi$. 
Let $L_\su$ be the limit set of $g_{\su}$. Then $L_\su$ is disjoint from
$g_{\su'}$ for all $\su' \in S(\su)$.
\end{lemma}

\proof
Suppose there exists a $\su'\in S(\su)$ and a $t>t_{\su'}$ with $g_{\su'}(t)
\in L_\su$. Since all $|\Re(E^{\circ k}(g_{\su}(t)))|>R$ for all $t>t_\su$ and
all $k\in\N$, all points $g_{\su}(t)$ have external address $\su$, and so has
the limit $g_{\su'}(t)$ (Theorem~\ref{fast}).

Now choose a potential $t' \in ]t_{\su},t[$.
By Proposition~\ref{control} (and possibly after finitely many iterations) we get
$|\Re(g_{\su'}(t))|-|\Re(g_{\su'}(t'))| \geq 4$. Then there is a potential $t''
>t_{\su}$ arbitrarily close to $t_{\su}$ such that $g_{\su}(t'')$ is
arbitrarily close to  $g_{\su'}(t)$. More precisely we assume that $t''<t'$ and
$|g_\su(t'')-g_{\su'}(t)|<1$, hence 
$ |\Re(g_{\su}(t''))|-|\Re(g_{\su'}(t'))|>3$. 
Since $g_\su(t'')$ and $g_{\su'}(t')$ both have external address $\su$, we
have $\left|\Im\left(E^{\circ k}(g_\su(t''))-E^{\circ
k}(g_{\su'}(t'))\right)\right|\leq 2\pi$ for all $k$. By
Lemma~\ref{LemExpoSeparation} this means $t''>t'$ and we get a contradiction.
Hence $L_\su \cap g_{\su'}(]t_{\su},\infty[) = \emptyset$ for all $\su' \in
S(\su)$.
\qed

\begin{theorem}[Escaping Points are Organized in Rays]
\label{organ} \lineclear
For every escaping point $w$  there exists a unique exponentially bounded
external address $\su$ and a unique potential $t \geq t_{\su}$ such that exactly
one of the following holds:
\begin{itemize}
\item 
either $t > t_{\su}$ and $w= g_{\su}(t)$,
\item  
or $t=t_{\su}$ and the dynamic ray $g_{\su}$  lands at $w$ such that
$w$ and the ray $g_{\su}$ escape uniformly,
\item 
or, if one of the singular values escapes: $v_i=g_{\su}(t)$ for
some $\su$ and $t \geq t_{\su}$, and the point $w$ maps to
$g_{\su}(t')$ with $t_{\su}\leq t'<t$ after finitely many iterations. 
\end{itemize}
\end{theorem}
\proof
Let $R$ be as defined in Lemma~\ref{LemEscapingConnected}. If the entire orbit
of $w$ lies in $\C \sm Y_{R+2}$, then let $\su$ be the external address of
$w$. By Lemma~\ref{LemEscapingConnected}, either $w=g_\su(t)$ for some
$t>t_\su$, or the ray $g_\su$ lands at $w$ with uniform escape.  In both
cases, the external address $\su$ is uniquely determined by the orbit of
$w$, and so is $t$ along dynamic rays because each ray is an injective curve.
To finish the uniqueness claim, we have to show that $g_\su(t)\neq
g_\su(t_\su)$ for all $t>t_\su$: this follows from
Lemma~\ref{LemLimitSetNotRay}.

If not the entire orbit of $w$ is in $\C \sm Y_R$, there exists a finite iterate
of $w$ whose orbit has that property. This iterate is either on a dynamic ray or
landing point of a ray.  By pulling back along the orbit of $w$ the claim is
proved for all $w$; if a singular value lies on a dynamic ray along this
pull-back, then an exception may occur as stated in the third case of the claim.
\qed

Now we want to show under which conditions a landing point escapes. We need slow
and fast external addresses, just like in \cite{Escaping}.

\begin{definition} [Slow and Fast External Addresses] \lineclear
We say that an external address $\su$ is {\em slow} if there are
$A,x>0$ and infinitely many $n$ for which $|s_{n+k}| \leq AF^{\circ
(k-1)}(x)$ for all $k \geq1$. Otherwise we call $\su$ {\em fast}.
\end{definition}

Note that every external address $\su$ with $t_\su>0$ is fast, but the converse
is not true: the two external addresses $1\,2\,1\,3\,1\,4\,1\,5\dots$ and
$1\,2\,1\,2\,3\,1\,2\,3\,4\,1\,2\,3\,4\,5\dots$ (with arbitrary entries $L$ or
$R$) are both unbounded with $t_\su=0$, but the first one is fast while the
second one is slow.

Now the following result holds in complete analogy to
\cite[Proposition~6.8]{Escaping}, so we omit the proof.

\begin{proposition} [Uniform Escape for Fast Addresses] 
\label{PropUniformEscape}\lineclear
An external address $\su$ is fast if and only if the ray $g_{\su}$
lands at an escaping point so that ray and landing point escape
uniformly.
\qedd
\end{proposition} 

\begin{corollary}[Uniform Escape of Ray and Landing Point]
\label{CorUniformEscape} \lineclear
If a dynamic ray lands at an escaping point, then the ray and its
landing point escape uniformly.
\end{corollary}
\proof
Let $g_\su$ be a dynamic ray and $w$ be its escaping landing point. If $\su$
is fast (in particular, if $t_{\su}>0$), then $g_\su$ lands at $w$ with
uniform escape by Proposition~\ref{PropUniformEscape}. Therefore, $\su$ is
slow and $t_\su=0$.

The point $w$ cannot be on any other dynamic ray $g_{\su'}$ by the same
argument as in \cite[Corollary~6.9]{Escaping}: every point $g_{\su'}(t)$ with
$t>t_{\su'}$ can be approximated arbitrarily closely by pieces of other rays
almost parallel to $g_{\su'}$ so that $g_{\su'}(t)$ is not accessible by any
ray. Therefore, by our classification, $w$ is the escaping landing point of
some other dynamic rays $g_{\su'}$, and thus $\su'$ is necessarily fast. This
means that the dynamic rays $g_{\su'}$ and $g_\su$ land together at $w$, so
for every $N\in\N$ there is a $k\in\N$ such that the first entries of
$\sigma^k(\su')$ and $\sigma^k(\su)$ differ by at least $N$, while
$g_{\sigma^k(\su')}$ and $g_{\sigma^k(\su)}$ land together at $E^{\circ k}(w)$.
But the $2\pi i n$-translates of the same rays must land at the $2\pi i
n$-translates of $E^{\circ k}(w)$, and this is a topological impossibility if
$N>1$.
\qed

As a result, we obtain the following classification of escaping points:

\begin{theorem}[Classification of Escaping Points]
\label{ThmClassEscaping} \lineclear
If no critical value escapes, then the set of escaping points is classified by
external addresses $\su\in\Sym$ so that 
\begin{itemize}
\item
for $\su$ with $t_\su<\infty$ so that $\su$ is slow, the associated set of
escaping points is the ray $g_\su(]t_\su,\infty[)$;
\item
for $\su$ with $t_\su<\infty$ so that $\su$ is fast, the associated set of
escaping points is $g_\su([t_\su,\infty[)$, i.e.\ the ray including its escaping
landing point with uniform escape;
\item
for $\su$ with $t_\su=\infty$, there is no associated set of escaping points.
\end{itemize}
The sets of escaping points to all $\su\in\Sym$ are disjoint, and their union
is the entire set of escaping points of $E$. 

Moreover, each path component of the set of escaping points is exactly the set
of escaping points associated to any particular exponentially bounded external
address.
\qedd
\end{theorem}
If one or both critical values escape, the necessary modifications are
straightforward.

\section{Epilogue: Hausdorff Dimension} 
\label{SecEpilogue}

\begin{theorem}[Dimension Paradox] 
\label{ThmDimension} \lineclear
The union of all dynamic rays has Hausdorff dimension 1,
while the set of the escaping landing points has Hausdorff dimension 2 and even
infinite planar Lebesgue-measure.
\end{theorem}
\sketch
McMullen~\cite{McMullen} has shown that the set of escaping points for every
$E_{a,b}$ has infinite planar Lebesgue measure. By our classification, every
escaping point is either on a dynamic ray or the landing point of a dynamic
ray. It thus suffices to prove that the union of all dynamic rays has
Hausdorff dimension $1$. We sketch a lemma of Karpi\'nska~\cite{Karpinska}
which generalizes to show that the set of rays has dimension less than $1+1/p$
for every $p>0$.

Given $p>0$ and $\xi>0$, let \[
P_{p,\xi}:=\{(x+iy)\in\C\colon |x|>\xi, |y|<|x|^{1/p}\} \,\,.
\]
We will show that $S_{p,\xi}:=\{z\in\C\colon E^{\circ k}(z)\in
P_{p,\xi} \mbox{ for all $k$}\}$ has Hausdorff dimension at most $1+1/p$ for
large $\xi$.

Let $Q\subset P_{p,\xi}$ be a square of side length $2\pi$ with boundaries
parallel to the coordinate axes, and with real parts in $[x,x+2\pi]$. If
$|x|>\xi$ is large, then $E(Q)$ is almost an annulus between radii $e^x$ and
$e^{2\pi}e^x$ (up to a factor $|a|$ or $|b|$ which will not matter). Then
$E(Q)\cap P_{p,\xi}$ is contained in two approximate rectangles of width
$e^x(e^{2\pi}-1)$ and height $2e^{x/p}e^{2\pi/p}$. Filling it with
squares of side length $2\pi$, we need approximately $Ae^{x(1+1/p)}$ squares
(where $A>0$ is some bounded factor). Pulling those back into the original
square $Q$, we obtain the same number of squares, with diameters approximately
$e^{-x}/2\pi$. Now
\[
\frac{\log\left(Ae^{x(1+1/p)}\right)}{-\log\left(e^{-x}/2\pi\right)}
=\frac{x(1+1/p)+\log(A)}{x+\log(2\pi)}=1+\frac 1 p+O(1/\xi)
\,\,.
\]
We can refine the tiling of $S_{p,\xi}$ by iterating this procedure: the
squares with side length $2\pi$ in $E(Q)\cap P_{p,\xi}$ can be replaced
themselves by smaller squares. In the limit, we obtain a covering of
$S_{p,\xi}$ with small tiles, and we obtain an upper bound of $1+1/p$ for the
Hausdorff dimension of $S_{p,\xi}$. 

In Proposition~\ref{control} we have shown that every point $z$ on a ray 
satisfies the parabola condition $|\Im(E^{\circ 
n}(z))|^p<|\Re(E^{ \circ n}(z))|$ for any $p>0$ for all but finitely many $n$.
Therefore, the union of all rays is contained in the countable union
\[
\bigcup_{n\geq 0} E^{-N}\left(S_{p,\xi}\right)
\]
and still has Hausdorff dimension at most $1+1/p$ for every $p>0$. Therefore,
the union of all rays has Hausdorff dimension $1$.
\qed

As mentioned in the introduction, the dynamics is particularly easy to
understand when both critical orbits are strictly preperiodic: see
\cite{Cosine}. In this case, {\em every} dynamic ray lands somewhere in $\C$,
and it can be read off from the external address whether the landing point is
periodic, recurrent, escaping, whether its orbits is dense in $\C$ etc..
Conversely, every point in $\C$ which is not on a dynamic ray is the landing
point of a ray. The set of dynamic rays still has Hausdorff dimension $1$, so
the set of landing points of this $1$-dimensional set of rays is the entire
complex plane minus this $1$-dimensional set of rays. In that case, almost
every $z\in\C$ (with respect to planar Lebesgue measure) is an escaping point.


\small
\begin{thebibliography}{McM}

\bibitem[DH]{DH}
Adrien Douady and John Hubbard:
{\em Etude dynamique des polyn\^omes complexes}.
Pr\'epublications math\'ematiques d'Orsay {\bf 2} (1984) and {\bf 4} (1985).

\bibitem[E]{E}
Alexandre Eremenko:
{\em On the iteration of entire functions}. In: {\em Dynamical 
systems and ergodic
theory}, Banach Center publications, Polish Scientific Publishers {\bf 23}
(1989), 339-345.

\bibitem[EL]{EL}
Alexandre Eremenko and Mikhail Lyubich:
{\em Dynamical properties of some classes of entire functions}.
Annales de l'Institut Fourier, Grenoble {\bf 42} 4 (1992),
989--1020.

\bibitem[K]{Karpinska}
Bogus{\l}awa Karpi\'nska: {\em Hausdorff dimension of the hairs
without endpoints for $\lambda \exp(z)$}. Comptes Rendus Acad.\
Sci.\ Paris {\bf 328} S\'erie~I (1999), 1039--1044.

\bibitem[McM]{McMullen}
Curt McMullen: {\em Area and Hausdorff dimension of Julia sets of
entire functions}. Transactions of the American Mathematical
Society {\bf 300} 1 (1987), 329--342.

\bibitem[Re]{Lasse}
Lasse Rempe: {\em Dynamics of exponential maps}. Ph.D. thesis,
Christian-Al\-brechts-Universit\"at Kiel (2003);
\verb*
http://e-diss.uni-kiel.de/diss_781.

\bibitem[Ro]{Guenter}
G\"unter Rottenfu{\ss}er: {\em Dynamic Rays for a family of entire
transcendental maps}. Diploma thesis, Technische Universit\"at M\"unchen
(2002).

\bibitem[S]{Cosine}
Dierk Schleicher: {\em The dynamical fine structure of iterated cosine maps}.
Manuscript, in preparation. 

\bibitem[SZ]{Escaping}
Dierk Schleicher and Johannes Zimmer: {\em Escaping points of 
exponential maps}.  Journal of the London Mathematical Society {\bf 67} (2003),
1--21.

\end{thebibliography}
\end{document}